% Logic Eprints
%Submitted 1911 Tue Nov 24, 1992 by: morass@math.mit.edu ()
%logic/friedman_sy/minimal_universes.tex
%

%MINIMAL UNIVERSES
%AUTHOR: Sy D. Friedman

\input amstex
\magnification=\magstep1
\documentstyle{amsppt}
\NoBlackBoxes
\hsize 6.0truein
\vsize 8.2truein
\hcorrection{.2truein}
\vcorrection{.2truein}
\loadbold       
\define\low{\langle 0^{\#}, C\rangle, i-\Sigma }
\define\lown{\langle 0^{\#}, C\rangle, i_n-\Sigma }

\define\ORD{\operatorname{ORD}}

\define\Lim{\operatorname{Lim}}

\define\Dom{\operatorname{Dom}}
\define\QQ{Q}

\document
\baselineskip=15pt
%\baselineskip=.25truein

\font\bigtenrm=cmr12 scaled\magstep2
\centerline
{\bigtenrm {MINIMAL UNIVERSES}}

\vskip20pt

\font\bigtenrm=cmr10 scaled\magstep2
\centerline
{\bigtenrm{Sy D. Friedman}\footnote"*"{Research supported by NSF contract \#
8903380. }} 

\centerline
{\bigtenrm {M.I.T.}}

\vskip20pt

\comment
lasteqno 1@0
\endcomment

We work in G\"odel-Bernays class theory. And we say that a structure 
$\langle M,A\rangle$ is a {\it model of $ZFC$ } if $M$ is a model of $ZFC$ 
and obeys replacement for formulas which are permitted to mention $A\subseteq
M$ as a unary predicate. An inner model $M$ is {\it minimal} if there is a
class $A$  such that $\langle M,A\rangle$ is amenable yet has no transitive
proper elementary submodel. $M$ is {\it strongly minimal on a club} if there
is a club $C$  such that $\langle M,C\rangle$ is amenable and $\alpha\in
C\longrightarrow$ $\langle V^M_\alpha,C\cap\alpha\rangle$ is not a model of
$ZFC.$  Strong minimality on a club implies minimality, but not conversely.
It is consistent for $L$  to be strongly minimal on $C=ORD$ and if
$0^{\#}$ exists, $L$ is not minimal yet $L[0^{\#}]$ may or may not be minimal.

If $M_1\subseteq M_2$ are inner models then $M_2$ is a {\it locally generic}
extension of $M_1$ if every $x\in M_2$ belongs to a set-generic extension
of $M_1.$  Our main result states that if $V$  is strongly minimal on a club
and $0^{\#}$ exists then  some inner model is both minimal and a locally
generic extension of $L.$  $V$  can always be made strongly minimal on
a club by forcing a strongly minimalizing club without adding sets (Theorem
1). Thus if $0^{\#}$ exists then there does exist an inner model which is both
minimal and a locally generic extension of $L,$ definable in a forcing
extension of $V$  that adds no sets. A special case is when $V=$ the minimal
model of $ZFC+0^{\#}$ exists, in which case there is an inner model which is
minimal and does not contain $0^{\#}.$ This answers a question of Mack
Stanley. 

\proclaim{Theorem 1} \ (Folklore) \ There is a class forcing to add a club
$C$  such that $\langle V,C\rangle$ is a model of $ZFC$ and $\alpha\in
C\longrightarrow\langle V_\alpha,C\cap\alpha\rangle$ is not a model of $ZFC.$
\endproclaim

\demo{Proof} Conditions are bounded closed sets $p$  such that $\alpha\in
p\longrightarrow\langle V_\alpha,p\cap\alpha\rangle$  is not a model of 
$ZFC.$  Conditions are ordered by end extension. To preserve $ZFC$ it's enough
to show that if $\langle D_i|i<\lambda \rangle$ is a $\Sigma_n$ definable
sequence of open dense classes then the intersection of the $D_i$'s is dense.
Given  a condition $p,$  first extend if necessary so that $p$  contains an
ordinal greater then $\lambda$ and the parameters defining 
$\langle D_i|i<\lambda \rangle$ and then build a canonical
$\Sigma_n$-elementary chain of models 
$\langle  V_{\alpha_i}|i<\lambda\rangle$  
and
extensions $p_i$ of $p$  in $V_{\alpha_{i+1}}-V_{\alpha_{i}}$ such that
$p_i$ meets $D_i.$ Then at limit stages $\bar\lambda\le\lambda,$
$p_{\bar\lambda}$  is a condition since $\alpha_{\bar\lambda }$ is
$V_{\alpha_{\bar\lambda }}$-definably singularized by
$\langle \alpha_i|i<\bar\lambda\rangle.$  \hfill{$\dashv$}
\enddemo

\proclaim{Theorem 2} \ Suppose $V$ is strongly minimal on a club and $0^{\#}$
exists. Then there is a minimal locally generic extension of $L.$  
\endproclaim

Theorem 2 is proved using backwards Easton forcing, where $0^{\#},C$ are used
to select the appropriate (minimal) almost generic extension, $C$  being a
strongly minimalizing club. We first describe the building blocks of this
backwards Easton iteration,  which are designed to produce ``generic stability
systems''. 

\vskip10pt

{\bf Definition.}  \  A {\it stability system} $p$  consists of a successor
ordinal $|p|=\alpha (p)+1$ and functions $f_{k}=f^p_k,k>0$ such that

(a) \  $\Dom f_1=\Lim\cap|p|,$ $f_1(\alpha)\le\alpha$  for $\alpha\in\Dom
f_1,$ 
$f_1(\alpha)=$ \underbar{lim} $\langle
f_1(\bar\alpha)|\bar\alpha\in\Lim\cap\alpha\rangle$  for
$\alpha\in\Lim^2\cap|p|.$  Define
$\alpha<_1\beta\Longleftrightarrow\alpha<\beta$  and $\alpha<\gamma\le\beta,$
$\gamma\in\Dom f_1\longrightarrow f_1(\gamma)\ge\alpha.$  Then $\alpha\in\Dom
f_1\longrightarrow f_1(\alpha)\le_1\alpha.$ 

(b) \ $\Dom f_{k+1}=\{\alpha<|p|\bigg|\alpha \, a<_k$-limit$\},$
$f_{k+1}(\alpha)\le\alpha$  for $\alpha\in\Dom\,  f_{k+1},$
$f_{k+1}(\alpha)=$ \underbar{lim} $\langle
f_{k+1}(\bar\alpha)|\bar\alpha<_k\alpha,$  $\bar\alpha\in\Dom f_{k+1}
\rangle$ for
$\alpha<|p|,$ $\alpha\in<_k-\lim^2.$ Define
$\alpha<_{k+1}\beta\longleftrightarrow\alpha<_k\beta$ and
$\alpha<\gamma\le_k\beta,$ $\gamma\in\Dom f_{k+1}\longrightarrow
f_{k+1}(\gamma)\ge\alpha.$  Then $\alpha\in\Dom f_{k+1}\longrightarrow
f_{k+1}(\alpha)\le_{k+1}\alpha.$  

Intuitively, $f_k(\alpha)$ represents the supremum of the ordinals which are
``$\Sigma_k$ stable in $\alpha$'', but only in a formal sense.

\vskip10pt

{\bf Definition.} \ Suppose $\kappa$ is regular, $\ell>0,$ 
$\gamma<\kappa.$  The
forcing ${\Cal{P}}(\kappa,\ell,\gamma)$ consists of all stability systems $p$
such that $\gamma\le_{\ell}\alpha(p)<\kappa.$  Extension of conditions is
defined by: $q\le p\longleftrightarrow f^q_k\supseteq f^p_k$ for all $k$ and
$\alpha(p)\le^q_{\ell-1}\alpha(q),$  where $\le^q_0=\le.$ We will see that
$\le$ is transitive. 

\vskip10pt

\proclaim{Lemma 1} For any stability system and $k;$ $\le^p_k$ is  a tree
ordering and $\alpha\leqslant\beta\le^p_k\gamma,$
$\alpha\le^p_{k+1}\gamma\longrightarrow\alpha\le^p_{k+1}\beta.$  Also
$\{\alpha|\alpha<^p_k\beta\}$ is closed in $\beta.$ 
\endproclaim

\demo{Proof} We first prove that $\le_k=\le^p_k$ is a tree ordering, by
induction on $k.$  For $k=0$ we define $\le_0=\le$ and the result is clear.
Suppose that the result holds for $k$  and we wish to show that $\le_{k+1}$ is
a tree ordering. Reflexivity is clear since we mean $\le_{k+1}$ to include
$=.$  Antisymmetry is clear since
$\alpha\le_{k+1}\beta\longrightarrow\alpha\le\beta.$ Suppose
$\alpha\le_{k+1}\beta\le_{k+1}\gamma$ and we want $\alpha\le_{k+1}\gamma.$
Since we have by definition $\alpha\le_k\beta\le_k\gamma$ by induction we know
$\alpha\le_k\gamma.$ Suppose $\alpha<\delta\le_k\gamma,$ $\delta\in\Dom
f_{k+1}.$  If $\delta>\beta$ then since $\beta\le_{k+1}\gamma$ we have
$f_{k+1}(\delta)\ge\beta\ge\alpha.$  If $\delta\le\beta$ then
$\delta\le_k\beta$ since $\le_k$ is a tree
ordering and both $\delta$ and $\beta$  are $\le_k\gamma.$  Since
$\alpha\le_{k+1}\beta$ we have $f_{k+1}(\delta)\ge\alpha.$  So we have shown
that $\le_{k+1}$ is transitive. Now suppose $\alpha\le\beta$ are both
$\le_{k+1}\gamma.$  By induction
 $\alpha\le_k\beta.$  If $\alpha<\delta\le_k\beta,$
$\delta\in\Dom f_{k+1}$ then $\delta\le_k\gamma$ since $\le_k$ is transitive
so $f_{k+1}(\delta)\ge\alpha$ since $\alpha\le_{k+1}\gamma.$ So
$\alpha\le_{k+1}\beta$  and we have shown that $\le_{k+1}$ is a tree ordering.

If $\alpha\leqslant\beta\le_k\gamma,$ $\alpha\le_{k+1}\gamma$ then
$\alpha\le_k\beta$ since $\le_k$ is a tree ordering. If
$\alpha<\delta\le_k\beta,$ $\delta\in\Dom f_{k+1}$  then $\delta\le_k\gamma$
since $\le_k$ is transitive so $f_{k+1}(\delta)\ge\alpha$ since
$\alpha\le_{k+1}\gamma.$  So $\alpha\le_{k+1}\beta.$ 

Finally we show that $\{\alpha|\alpha<_k\beta\}$ is closed in $\beta,$ by
induction on $k.$  This is clear for $k=0.$  Suppose it holds for $k$  and
$\bar\alpha$ is a limit of $\{\alpha|\alpha<_{k+1}\beta\},\bar\alpha<\beta.$ 
Then $\bar\alpha<_k\beta$ by induction. Suppose
$\bar\alpha<\gamma\le_k\beta,\gamma\in\Dom f_{k+1}.$  Then
$f_{k+1}(\gamma)\ge\alpha$  for all $\alpha<_{k+1}\beta,$ $\alpha<\gamma.$  In
particular this holds when $\alpha<\bar\alpha$  so
$f_{k+1}(\gamma)\ge\bar\alpha$ since $\bar\alpha$ is a limit of such $\alpha.$
So $\bar\alpha<_{k+1}\beta.$  \hfill{$\dashv$}

\enddemo

\proclaim{Lemma 2} Let $f_k,$ $<_k$ arise from a stability system. If
$\alpha\in\Dom f_k,$ $f_k(\alpha)<\alpha$ then $f_k(\alpha)=$ largest
$\bar\alpha<_k\alpha.$ If $\alpha\in<_k-\lim^2,$ $f_{k+1}(\alpha)=\alpha$ then
$\{\bar\alpha|\bar\alpha<_{k+1}\alpha\}$ is unbounded in $\alpha.$ 
\endproclaim

\demo{Proof} Suppose $\alpha\in\Dom f_k,$ $f_k(\alpha)<\alpha.$ We know that
$f_k(\alpha)<_k\alpha,$ by definition of condition. If
$f_k(\alpha)<\beta<\alpha$ then $\beta\nless_k\alpha$ since
$f_k(\alpha)\ngeq\beta.$  So $f_k(\alpha)=$ largest $\bar\alpha<_k\alpha.$ 
Suppose $\alpha\in<_k-\Lim^2,$ $f_{k+1}(\alpha)=\alpha.$ Then 
$f_{k+1}(\alpha)=\alpha=\underline{\text{lim}}\{f_{k+1}(\bar\alpha)|
\bar\alpha<_k\alpha,$
$\bar\alpha\in<_k-\Lim\}.$  So
$\alpha=\lim\{f_{k+1}(\bar\alpha)|\bar\alpha<_k\alpha,\bar\alpha\in<_k-\Lim\}.$  For
any $\alpha_0<\alpha$ there must be $\bar\alpha_0<_k\alpha,$
$\bar\alpha_0\in<_k-\lim$ such that $f_{k+1}(\bar\alpha)\ge
f_{k+1}(\bar\alpha_0)\ge\alpha_0$  for all
$\bar\alpha<_k\alpha,\bar\alpha\in<_k-\lim,$  $\bar\alpha\ge\bar\alpha_0.$
(For, we need only first choose $\bar\alpha'_0$ to guarantee
$f_{k+1}(\bar\alpha)\ge\alpha_0$  for all $\bar\alpha$ beyond $\bar\alpha'_0$
and then minimize $f_{k+1}(\bar\alpha_0')$ to get $\bar\alpha_0.)$ Then
$f_{k+1}(\bar\alpha_0)<_k$ $\alpha$ since either
$f_{k+1}(\bar\alpha_0)<_{k+1}\bar\alpha_0<_k\alpha$ or
$f_{k+1}(\bar\alpha_0)=\bar\alpha_0<_k\alpha.$  Suppose
$f_{k+1}(\bar\alpha_0)<\beta\le_k\alpha,$ $\beta\in\Dom f_{k+1}.$ If
$\beta\le\bar\alpha_0$ then $f_{k+1}(\beta)\ge f_{k+1}(\bar\alpha_0)$ since
$f_{k+1}(\bar\alpha_0)<_{k+1}\bar\alpha_0.$  If $\bar\alpha_0\le\beta$ then
$\bar\alpha_0\le_k\beta$ and $f_{k+1}(\beta)\ge f_{k+1}(\bar\alpha_0)$  by
choice of $\bar\alpha_0.$  So $f_{k+1}(\bar\alpha_0)<_{k+1}\alpha$  and
$f_{k+1}(\bar\alpha_0)\ge\alpha_0.$  So
$\{\bar\alpha|\bar\alpha<_{k+1}\alpha\}$  is unbdd in $\alpha.$ \hfill{$\dashv$}
\enddemo

\proclaim{Lemma 3} Suppose $r\le q\le p$ in ${\Cal{P}}(\kappa,\ell,\gamma).$
Then $r\le p.$  
\endproclaim

\demo{Proof} We need to check that $\alpha(p)\le^r_{\ell-1}\alpha(r).$ But
$\alpha(p)\le^q_{\ell-1}\alpha(q)$  and so $\alpha(p)\le^r_{\ell-1}\alpha(q),$
and $\alpha(q)\le^r_{\ell-1}\alpha(r).$  So the result follows from Lemma 1.
\hfill{$\dashv$} 
\enddemo

\proclaim{Lemma 4} Suppose $p\in {\Cal{P}}(\kappa,\ell,\gamma)$  and
$\alpha(p)\le\alpha<\kappa.$ Then there exists $q\le p,$ $\alpha(q)=\alpha.$ 
\endproclaim

\demo{Proof} For limit $\lambda\in(\alpha(p),\alpha ]$ define
$f^q_k(\lambda)=\lambda.$  It is routine to verify that the resulting $q$ is a
condition and extends $p.$  \hfill{$\dashv$}
\enddemo

\proclaim{Lemma 5} Suppose $p_0\ge p_1\ge\cdots$ is a sequence of conditions
in 
${\Cal{P}}(\kappa,\ell,\gamma)$  of length $<\kappa.$  Then there is $p\le$
each $p_i,$ $\alpha(p)=\bigcup\limits_{i}\alpha(p_i).$ 
\endproclaim

\demo{Proof} Assume that the $p_i$'s are distinct. Let
$\alpha=\bigcup\limits_{i}\alpha(p_i).$  We must define $f^p_k(\alpha).$  We
do so by induction
 on $k>0.$ If $\alpha\notin\Lim,$ $f^p_1(\alpha)$ is undefined.
If $\alpha\in\Lim^2,$ let $f^p_1(\alpha)=$ \underbar{lim} $\langle
f^{p_{i}}_1(\bar\alpha)|\bar\alpha\le\alpha(p_i),$ $\bar\alpha\text{
limit}\rangle.$  If $\alpha\in\Lim-\Lim^2$ then let $f^p_1(\alpha)=\alpha.$
Assuming $f^p_k(\alpha)$ is defined (and
$f^p_k\restriction\alpha=\bigcup\limits_i f^{p_{i}}_k)$ it makes sense to ask
if $\alpha\in<^p_k-\lim.$ If not then $f^p_{k+1}(\alpha)$ is undefined. If
$\alpha\in<^p_k-\lim^2$ then set $f^p_{k+1}(\alpha)=$ \underbar{lim} $\langle
f^p_{k+1}(\bar\alpha)|\bar\alpha<^p_k\alpha,\bar\alpha \text{ a
}<^p_k-\text{limit}\rangle.$ If $\alpha\in<^p_k-\lim-<^p_k-\lim^2$ then set
$f^p_{k+1}(\alpha )=\alpha.$

Now we show that $\alpha(p_i)<^p_{\ell-1}\alpha,$ defined in terms of the
above $f^p_k$'s. Suppose $\alpha(p_i)<^p_k\alpha$ for all $i,$  where
$k<\ell-1$ and we want $\alpha(p_i)<^p_{k+1}\alpha.$ Suppose
$\alpha(p_i)<\beta\le^p_k\alpha,$ $\beta\in\Dom f^p_{k+1}.$  If $\beta<\alpha$
then we can choose $j$  so that $\beta\le^{p_{j}}_k\alpha(p_j)$  and then
$f^{p_{j}}_{k+1}(\beta)\ge\alpha(p_i)$  since $p_j\le p_i.$ If $\beta=\alpha$
then $f^p_{k+1}(\beta)<\alpha(p_i)$  can only result if
$f^p_{k+1}(\bar\alpha)<\alpha(p_i)$ for some $\bar\alpha<^p_k\alpha,$
$\bar\alpha>\alpha(p)$ but then $\bar\alpha<^{p_{j}}_k\alpha(p_j)$ for large
$j,$   contradicting $p_j\le p_i.$  So $\alpha(p_i)\le^p_{\ell-1}\alpha$  for
all $i.$  

Now it is easy to verify that $p$ is a condition extending each $p_i,$ since
any violation created by $\alpha=\alpha(p)$ would imply a violation at some
$\alpha(p_i)<^p_{\ell-1}\alpha(p).$  \hfill{$\dashv$}
\enddemo

Now we describe the backwards Easton iteration used to create our minimal
inner model. ${\Cal{P}}$ is the iteration with Easton supports over $L$  where
${\Cal{P}}_0=$ the trivial forcing, ${\Cal{P}}_\lambda=$ inverse limit at
singular $\lambda,$  direct limit at regular $\lambda,$  
${\Cal{P}}_{\kappa +1}={\Cal{P}}_\kappa *\dot\QQ_\kappa$ where
$\dot\QQ_\kappa$ is a term for the trivial forcing unless $\kappa$ is regular.
For regular $\kappa,$ $\dot\QQ_\kappa$ is a term for the following forcing in
$L[G_\kappa ],$ $G_\kappa$ denoting the ${\Cal{P}}_\kappa$-generic: choose a
pair $(\ell_\kappa,\gamma_k)$ with $\ell_\kappa>0,$ $\gamma_\kappa<\kappa$
and apply the forcing ${\Cal{P}}(\kappa,\ell_\kappa,\gamma_\kappa).$  Now
${\Cal{P}}_\kappa\Vdash\dot\QQ_\kappa$ is $<\kappa$-closed and has cardinality
$\kappa,$ so ${\Cal{P}}$ preserves cofinalities.

Our goal is to build $G=\langle G_\alpha|\alpha\in\text{ ORD }\rangle$ so that
$G_\alpha$ is ${\Cal{P}}_{\alpha }$-generic over $L$  and to select ordinals
$\alpha_i\in[i,i^*),$ $i<i^*$ adjacent Silver indiscernibles such that
(writing $G_{\alpha+1}=G_\alpha *g_\alpha):$ 

\vskip10pt

1. \ $i<j$ in $I=$ Silver indiscernibles, $p\in g_i,$ $q\in g_j,$
$\alpha(q)\ge i\longrightarrow f^p_k\subseteq f^q_k$  for all $k.$  let
$f_k=\bigcup\{f^p_k|p\in g_i$  for some $i\in I\}.$ 

\vskip10pt

2. \ For $i\in I,\ell_i=$ least $\ell$ such that the $\langle
0^{\#},C\rangle,i-\Sigma_\ell$ stables in $C$  are bounded in $i,$  where $C=$
the given strongly minimalizing club for $V.$ $(\alpha$ is
$B,\beta-\Sigma_\ell$ stable if $\langle L_\alpha [B], B\cap\alpha\rangle$ is
a $\Sigma_\ell$-elementary submodel of $\langle L_\beta
[B],B\cap\beta\rangle).$  Also $\gamma_i=\alpha_j$ where $j=\bigcup\{\langle
0^{\#},C\rangle,i-\Sigma_{\ell_{i}}$ stables in $C\}\ge 0.$  (By convention,
$\alpha_0=0.$)  

\vskip10pt

3. \ For $i\in I,$ $f_k(i)=i$ if $k<\ell_i$  and $f_{\ell_{i}}(i)=\gamma_i.$ 

\vskip10pt

4. \ For $i\in I, i\le_{\ell_{i}}\alpha_i$ (where $\le_k$ is defined from the
$f_k$'s) and $\alpha_i\in\Dom f_{\ell_{i}+1},
f_{\ell_{i}+1}(\alpha_i)=\alpha_j$  where $j=\bigcup\{\langle
0^{\#},C\rangle,i-\Sigma_{\ell_{i}+1}$ stables in $C\}.$

Suppose that the $f_k$'s have been constructed to obey 1--4 above and we now
prove Theorem 2. The desired minimal, locally generic extension of $L$ is
$L[\langle G_\alpha|\alpha\in\ORD\rangle ],$ witnessed by the amenable class
$\langle f_k|k\in\omega\rangle.$ The reason for minimality is roughly as
follows: there are unboundedly many $\alpha<_k\infty$ (defined from the
$f_k$'s) yet no ordinal $\alpha$  is $<_k\infty$ for all $k$  simultaneously.
More precisely:

\proclaim{Lemma 6} let $i<j$ be indiscernibles, $i\langle
0^{\#},C\rangle,j-\Sigma_k$ stable, $i\in C$ and ($j$ a limit of $\langle
0^{\#},C\rangle,j-\Sigma _{k-2}$ stables or $k\le 2).$  Then
$i\le_k\alpha_i<_k j,$  where $\le_k$  is defined from the $f_k$'s.
\endproclaim 

\demo{Proof}  By induction on $k$  and for fixed $k$  by induction on $j.$  By
property 4, $i\le_{\ell_i}\alpha_i$ and clearly $\ell_i\ge k$ since $i$ is a
limit of $\langle 0^{\#},C\rangle,i-\Sigma_{k-1}$ stables $(k>1)$ and so
$\ell_i\ge k$ follows by property 2. So we only need to check that
$\alpha_i<_k j.$  

Suppose $k=1.$  If $\alpha_i\nless_1 j$ then choose $\alpha_i<\beta\le j$ so
that $f_1(\beta)<\alpha_i.$  Then $\beta<j$ as properties 2, 3 imply that
$f_1(j)\ge\alpha_i,$ since $i$ is $\langle 0^{\#},C\rangle,j-\Sigma_1$ stable
and belongs to $C.$  There are no indiscernibles between $\beta$  and $j$  as
otherwise we can apply induction on $j.$  So $\bar j\le\beta< j$ where $\bar
j=I$-predecessor to $j.$  In fact $\bar j<\beta<j$ since otherwise
$\alpha_i<\beta=\bar j$ and again 2, 3 imply that $f_1(\bar j)\ge\alpha_i.$ If
$f_1(j)<\beta$ then since $f_1(j)<_1j$ we have $f_1(j)\le f_1(\beta)$ and
hence $f_1(\beta)\ge f_1(j)\ge\alpha_i,$ contrary to assumption. So
$f_1(j)\ge\beta>\bar j$ and by 3, $f_1(j)=\alpha_{\bar j},\bar j$ is $\langle
0^{\#},C\rangle,j-\Sigma_1$ stable and belongs to $C.$  And $\bar
j<\beta\le\alpha_{\bar j}.$  But $\bar j\le_1\alpha_{\bar j}$  so
$f_1(\beta)\ge\bar j$  and $\bar j>\alpha_i$ since $\beta\le\alpha_{\bar j}$
and $\beta>\alpha_i.$  So $f_1(\beta)>\alpha_i,$ contradicting our assumption.

Suppose the lemma holds for $k$  and we prove it for $k+1.$  If
$\alpha_i\nless_{k+1}j$ then choose $\alpha_i<\beta\le_k j$ so that
$\beta\in\Dom f_{k+1},$ $f_{k+1}(\beta)<\alpha_i.$  By 2, 3 we have $\beta<_k
j.$   By induction on $j$  there can be no $\bar j$  such that $\beta\le\bar j
$  and $\bar j$ is $\langle 0^{\#},C\rangle,j-\Sigma_k$ stable and in
$C,$ as  otherwise induction on $k$  implies $\bar j\le_k\alpha_{\bar j}<_kj$
so $\bar j<_k j$ and $\beta\le_k\bar j$ by Lemma 1. By 2, 3 $f_k(j)$ is
defined and equal to $\alpha_{\bar j}$ where $\bar j=\bigcup\{\langle
0^{\#},C\rangle,j-\Sigma_k$ stables in $C\}$ and since $\beta<_k j$ we have
$\beta\le\alpha_{\bar j}.$  And $\bar j<\beta\le\alpha_{\bar j}$ since there
is no $\bar j$  such that $\beta\le\bar j$  and $\bar j$ is $\langle
0^{\#},C\rangle,j-\Sigma_k$ stable, $\bar j\in C.$ Now $\alpha_{\bar
j}=f_k(j)<_k j$ so since $\beta<_k j$ we have $\beta\le_k\alpha_{\bar j}$  by
Lemma 1. Now let $\ell=\ell_{\bar j}.$ Clearly $\ell\ge k$ since $\bar j$ is
$\langle 0^{\#},C\rangle,j-\Sigma_k$ stable, $\bar j\in C$  and hence $\bar
j=\bigcup\{
\langle 0^{\#},C\rangle,\bar j-\Sigma_{k-1}$  stables in $C\}.$  But if
$\ell>k$ then by 4, $\bar j<_{k+1}\alpha_{\bar j}$ and this contradicts
$f_{k+1}(\beta)<\alpha_i<\bar j.$  So $\ell=k,$ $\bar j<_k$ $\alpha_{\bar j}$
and by 4, $f_{k+1}(\alpha_{\bar j})$ is defined and equal to
$\alpha_{\overset=\to j}$ 
where $\overset=\to j=\bigcup
\{\langle 0^{\#},C\rangle,\bar j-\Sigma_{k+1}$ stables in $C\}.$ 
But $i$ is $\langle
0^{\#},C\rangle,j-\Sigma_{k+1}$  stable, $i\in C$ and $\bar j$ is $\langle
0^{\#},C\rangle,j-\Sigma_k$ stable, $\bar j\in C$  and $i<\bar j$ so $i$ is
$\langle 0^{\#},C\rangle,\bar j-\Sigma_{k+1}$ stable and we get
$i\le\overset=\to j.$  Thus $\alpha_i\le f_{k+1}(\alpha_{\bar j})<\beta.$
This contradicts $\beta\le_k\,\alpha_{\bar j}, f_{k+1}(\beta)<\alpha_i$  since
$f_{k+1}(\alpha_{\bar j})<_{k+1}\alpha_{\bar j}.$  \hfill{$\dashv$}

\enddemo

\proclaim{Corollary 7} Define
$\alpha<_{k}\infty\longleftrightarrow\alpha<_k\beta$  for cofinally many
$\beta<_{k-1}\infty$  (where $\beta<_0\infty$ is vacuous). Then for each $k$
there are cofinally many $\alpha<_k\infty.$ 
\endproclaim

\demo{Proof} By Lemma 6 if $i$  is $\langle 0^{\#},C\rangle-\Sigma_k$ stable
and belongs to $C$  then $i<_k\infty$ (by induction on $k).$  The class of all
such $i$  is cofinal in the ordinals. \hfill{$\dashv$}
\enddemo

\proclaim{Corollary 8} No $\alpha$ is $<_k\infty$ for all $k.$ 
\endproclaim

\demo{Proof} Choose $k$ large enough so that $i=$ the least 
$\langle 0^{\#},C\rangle-\Sigma_k$ stable  in $C$  is larger than $\alpha.$
Then $i<_k\infty$ is $f_k(i)=\alpha_0=0.$  So $\alpha \nless_k\infty.$ 
\hfill{$\dashv$}
\enddemo

Thus we have minimality, since if 
$\langle L[\langle G_\alpha|\alpha\in\ORD\rangle ],$ $\langle
f_k|k\in\omega\rangle\rangle$ had a transitive elementary submodel, by
Corollary 7 its height $\alpha$  would be $<_k\infty$ for each $k,$ in
contradiction to Corollary 8.

It remains to construct the $G_\alpha$'s so as to obey 1--4.

\vskip10pt

\flushpar
{\bf The Construction.} \ We build $G_{i+1}=G_i*g_i$ by induction on $i\in I.$
When defining $G_{i^*+1}$ we also specify $\alpha_i\in[i,i^*),$  where $i^*$
denotes the $I$-successor to $i.$ 

\vskip10pt

\flushpar
$\bold{G_{i_0+1}, i_0=\min I}$ \  Choose  $G_{i_{0}}$ to be the
$L[0^{\#}]$-least generic for ${\Cal{P}}_{i_{0}},$ using the countability
of $i_0.$  Set $\ell_{i_{0}}=1,\gamma_{i_{0}}=0=\alpha_0$ and choose
$g_{i_{0}}$ to be the $L[0^{\#}]$-least generic for ${\Cal{P}}(i_0,1,0)$ as
defined in $L[G_{i_{0}}].$ 

\vskip10pt

\flushpar
$\bold{G_{i^*+1},i\in I}$  \  First choose $G_{i^{*}}$ to be the
$L[0^{\#},C]$-least generic for ${\Cal{P}}_{i^*}$ extending $G_{i+1},$ using
the $\le i$-closure of ${\Cal{P}}_{i+1,i^*}$ in $L[G_{i+1}],$ where
${\Cal{P}}_{i^*}={\Cal{P}}_{i+1}*{\Cal{P}}_{i+1,i^*}.$  (Note that the dense
sets in ${\Cal{P}}_{i+1,i^*}$  can be grouped into countably many collections
of size $i,$  enabling an easy construction of a generic.) The key step
involves the choice of $g_{i^*}.$ 

Choose $n\ge 0$ so that the ordertype of the $\low_\ell$ stables in $C$ is
$\lambda+n,$ $\lambda$ limit or $0,$  where $\ell=\ell_i=$ least $\ell$ such
that the $\low_\ell$ stables in $C$ are bounded in $i.$  Let
$\overline{Q}_{i^*}$ be the forcing ${\Cal{P}}(i^*,\ell+1,\alpha_j)$ as
defined in $L[G_{i^*}],$ where $j=\bigcup
\{\low_{\ell+1}$ stables in $C\}\ge 0.$ 
Let $p_0\in\overline{Q}_{i^*}$ be the ``condition'' defined by
$\alpha(p_0)=i,$ $f^{p_{0}}_k\restriction i= \bigcup\{f^p_k|p\in g_i\},$
$f^{p_{0}}_k(i)=i$ if $k<\ell,$ $f^{p_{0}}_\ell(i)=\gamma_i=\alpha_{j'}$ where
$j'=\bigcup\{\low_{\ell}$ stables in $C\}\ge 0.$  {\it We will verify} later
that $p_0$ is indeed a condition in $\overline{Q}_{i^*}.$  Choose $p_1\le p_0$
in $\overline{Q}_{i^*}$ 
to meet all dense $\Delta$ in $L[G_{i^{*}}]$ definable in
$L[G_{i^{*}}]$  from $G_{i^{*}}$  and parameters in
$(i+1)\cup\{j_1\cdots j_n\}$ where $j_1\cdots j_n$ are the first $n$
indiscernibles $\ge i^*.$  Also arrange that $\alpha(p_1)$ is a
$<^{p_{1}}_\ell-$limit, $f^{p_{1}}_{\ell+1}(\alpha(p_1))=\alpha_j.$ Now
set $\alpha_i=\alpha(p_1)$ and $\ell_{i^{*}}=1,$ $\gamma_{i^{*}}=\alpha_i.$
Choose $g_{i^{*}}$ to be generic for $Q_{i^*}={\Cal{P}}(i^*,1,\alpha_i)$ over
$L[G_{i^{*}}],$ extending the condition $p_1\in Q_{i^{*}}.$ 

\vskip10pt

\flushpar
$\bold{G_{i+1},i\in\Lim I}$ \  $G_i=\bigcup\{G_j|j\in I\cup i\}.$ Let
$Q_i={\Cal{P}}(i,\ell_i,\gamma_i)$ in $L[G_i]$ where $\ell_i$ is the least
$\ell$ such that the $\low_\ell$ stables in $C$ are bounded in $i$ and
$\gamma_i=\alpha_j$ where $j=\bigcup\{\low_{\ell_{i}}$ stables in $C\}.$ Let
$f_k\restriction i=\bigcup\{f^p_k|p\in g_j$ for some $j\in I\cap i\}.$  Then
$g_i=Q_i$-generic determined by $\langle f_k\restriction i 
|k\in\omega\rangle.$ 
{\it We will verify} later that $\langle f_k\restriction i|k\in\omega\rangle$
does indeed determine a $Q_i$-generic over $L[G_i].$ 

\proclaim{Lemma 9} Assume that the verifications claimed in the construction
can be carried out. Then 1--4 hold. 
\endproclaim 

\demo{Proof} Everything is clear, with the possible exception of the first
statement in 4: \  $i\le_{\ell_{i}}\alpha_i.$  But note that in the
construction of $G_{i^*+1,}$ $i=\alpha(p_0),$ $\alpha_i=\alpha(p_i)$ where
$p_1\le p_0$ in ${\Cal{P}}(i^*,\ell_i+1,\alpha_j)$ for some $j,$  so we're
done by the definition of extension of conditions. \hfill{$\dashv$}
\enddemo

\proclaim{Lemma 10} By induction on $i\in I:$ 

(a) \ $g_i$ is well-defined and $Q_i$-generic, where
$Q_i={\Cal{P}}(i,\ell_i,\gamma_i)$ as interpreted by $G_i.$

(b) \ Define $p_0$ by: \  $\alpha(p_0)=i,$ $f^{p_{0}}_k\restriction
i=\bigcup\{f^p_k|p\in g_i\},$ $f^{p_{0}}_k(i)=i$ if $k<\ell_i,$
$f^{p_{0}}_{\ell_{i}}(i)=\gamma_i.$  Then $p_0$ is a stability system.

(c)  \  Lemma 6 holds for indiscernibles $\le i.$ 

(d)  \  If $p_0$ is defined as in (b) then $p_0\in\overline{Q_{i^*}},$ as
defined in the construction.
\endproclaim

\demo{Proof} (a) \  This follows by induction unless $\ell_i>1$  and there is
a final segment $i_0<i_1<\cdots$ of the $\low_{\ell_i-1}$ stables $C$  of
ordertype $\omega.$  We may also assume that $i_0$ is big enough so that
$j=\bigcup
\{\low_{\ell_i}$ stables in $C\}=\bigcup
\{\langle 0^{\#},C\rangle,i_n-\Sigma_{\ell_i}$ 
stables in $C\}$ for all $n.$ Note that $\ell_{i_{n}}=\ell_i-1$ and the
ordertype of the $\lown_{\ell_{i_{n}}}$ stables in $C=\lambda+n'$ with $n\le
n'<\omega,$ $\lambda$ limit or $0.$  By construction, $p^n_1$ meets $\Delta$
in $L[G_{i^{*}_{n}}]$ defined from $G_{i^{*}_{n}}$ and parameters in $(i_n+1)
\bigcup$  (least $n$  indiscernibles $\ge i^*_{n})$ where $\Delta$ is dense on
${\Cal{P}}(i^*_n,\ell_i,\gamma_i)$  and $p^n_1$ is the condition in
$g_{i^{*}_{n}}$ with $\alpha(p^n_1)=\alpha_{i_{n}}.$  By an inductive use of
(c), $i_n<^{p_{0}}_{\ell_{i}-1}\alpha_{i_{n}}<^{p_{0}}_{\ell_i-1}i_m$ for all
$m>n,$  where $p_0$ is defined as in (b). So
$\alpha_{i_{0}}<^{p_{0}}_{\ell_{i}-1}\alpha_{i_{1}}<^{p_{0}}_{\ell_i-1}\cdots$
and $p^0_1\ge p^1_1\ge p^2_1\ge\cdots$ in ${\Cal{P}}(i,\ell_i,\gamma_i)$
determine the generic $g_i$  containing the $p^n_1$'s.

(b) \  The genericity of $g_i$ established in (a) implies that $i$ is a
$<^{p_{0}}_{\ell_i-1}-\lim^2$ and $f^{p_{0}}_k(i)$ is determined correctly by
$f^{p_{0}}_k(\gamma),$ $\gamma<^{p_{0}}_{k-1}i,$ $\gamma$ a
$<^{p_{0}}_{k-1}-\lim,$ for $k\le\ell_i.$  So $p_0$ obeys the requirements for
a stability system.

(c) \ The proof of Lemma 6 for indiscernibles $\le i$ only used the facts that
$p_0$  is a stability system and 1--4 hold $\le i.$  So we are done by Lemma
9 through $i.$ 

(d)  \  We must verify that $\alpha_j<^{p_{0}}_{\ell+1} i$ in the definition
of $g_{i^{*}}.$  This follows from (c).  \hfill{$\dashv$}
\enddemo

\enddocument